\documentclass[11pt]{article}
\usepackage{exscale,relsize}
\usepackage{amsmath}
\usepackage{amsfonts}
\usepackage{amssymb}
\usepackage{calc}
\usepackage{theorem}
\usepackage{pifont}      
\usepackage{xcolor}
\newcommand{\F}{\ensuremath{{F}}}

\newcommand{\emp}{\ensuremath{\varnothing}}

\newcommand{\scal}[2]{\langle{{#1},{#2}}\rangle}

\newcommand{\pscal}{\ensuremath{\scal{\cdot}{\cdot}}}

\newcommand{\RR}{\ensuremath{\mathbb R}}

\newcommand{\RP}{\ensuremath{\left[0,+\infty\right[}}
\newcommand{\RPP}{\ensuremath{\,\left]0,+\infty\right[}}
\newcommand{\RX}{\ensuremath{\,\left]-\infty,+\infty\right]}}

\newcommand{\NN}{\ensuremath{\mathbb N}}
\newcommand{\nnn}{\ensuremath{{n \in \NN}}}
\newcommand{\thalb}{\ensuremath{\tfrac{1}{2}}}

\newcommand{\menge}[2]{\big\{{#1} \mid {#2}\big\}}

\newcommand{\dom}{\ensuremath{\operatorname{dom}}}

\newcommand{\gra}{\ensuremath{\operatorname{gra}}}

\newcommand{\inte}{\ensuremath{\operatorname{int}}}

\newcommand{\bd}{\ensuremath{\operatorname{bdry}}}
\newcommand{\ran}{\ensuremath{\operatorname{ran}}}

\newcommand{\Fix}{\ensuremath{\operatorname{Fix}}}

\newcommand{\Id}{\ensuremath{\operatorname{Id}}}

\newcommand{\pinf}{\ensuremath{+\infty}}

\renewcommand{\phi}{\ensuremath{\varphi}}

\newcommand{\To}{\ensuremath{\rightrightarrows}}

\newtheorem{theorem}{Theorem}[section]
\newtheorem{lemma}[theorem]{Lemma}
\newtheorem{fact}[theorem]{Fact}
\newtheorem{corollary}[theorem]{Corollary}
\newtheorem{proposition}[theorem]{Proposition}
\newtheorem{definition}[theorem]{Definition}

\theoremstyle{plain}{\theorembodyfont{\rmfamily}
}
\theoremstyle{plain}{\theorembodyfont{\rmfamily}
}
\theoremstyle{plain}{\theorembodyfont{\rmfamily}
}
\theoremstyle{plain}{\theorembodyfont{\rmfamily}
\newtheorem{example}[theorem]{Example}}
\theoremstyle{plain}{\theorembodyfont{\rmfamily}
\newtheorem{remark}[theorem]{Remark}}
\theoremstyle{plain}{\theorembodyfont{\rmfamily}
}

\newcommand{\boxedeqn}[1]{%
    \[\fbox{%
        \addtolength{\linewidth}{-2\fboxsep}%
        \addtolength{\linewidth}{-2\fboxrule}%
        \begin{minipage}{\linewidth}%
        \begin{equation}#1\end{equation}%
        \end{minipage}%
      }\]%
  }

\begin{document}

\title{{\sffamily General Resolvents for Monotone Operators:\\
Characterization and Extension}}

\author{
Heinz H.\ Bauschke\thanks{Mathematics, Irving K.\ Barber School,
UBC Okanagan, Kelowna, British Columbia V1V 1V7, Canada. E-mail:
\texttt{heinz.bauschke@ubc.ca}.},~  Xianfu
Wang\thanks{Mathematics, Irving K.\ Barber School, UBC Okanagan,
Kelowna, British Columbia V1V 1V7, Canada. E-mail:
\texttt{shawn.wang@ubc.ca}.},~ and ~ Liangjin\
Yao\thanks{Mathematics, Irving K.\ Barber School, UBC Okanagan,
Kelowna, British Columbia V1V 1V7, Canada.
E-mail:  \texttt{ljinyao@interchange.ubc.ca}.} }

\date{
October 20, 2008 evening (\texttt{081020n.tex}) %
}

\maketitle


\begin{abstract} \noindent
Monotone operators, especially in the form of subdifferential
operators, are of basic importance in optimization.
It is well known since Minty, Rockafellar, and Bertsekas-Eckstein
that in Hilbert space, monotone operators can be understood
and analyzed from the alternative
viewpoint of firmly nonexpansive mappings,
which were found to be precisely the resolvents of monotone operators.
For example, the  proximal mappings in the sense of Moreau are
precisely the resolvents of subdifferential operators.
More general notions of ``resolvent'', ``proximal mapping''
and ``firmly nonexpansive'' have been studied.
One important class, popularized chiefly by Alber and by Kohsaka and 
Takahashi, is based on the normalized duality mapping.
Furthermore, Censor and Lent pioneered
the use of the gradient of a well behaved convex functions 
in a Bregman-distance based framework. 
It is known that resolvents are firmly nonexpansive,
but the converse has been an open problem for the latter framework.

In this note, we build on the very recent characterization
of maximal monotonicity due to Mart\'{\i}nez-Legaz
to provide a framework for studying
resolvents in which firmly nonexpansive mappings are always resolvents.
This framework includes classical resolvents, resolvents
based on the normalized duality mapping, resolvents based on
Bregman distances, and even resolvents based on (nonsymmetric)
rotators. As a by-product of recent work on the proximal average,
we obtain a constructive Kirszbraun-Valentine extension
result for generalized firmly nonexpansive mappings.
Several examples illustrate our results.
\end{abstract}

{\small 



}

\section{Introduction}

Throughout this paper, we assume that
$X$ is a real reflexive Banach space, with continuous dual space 
$X^*$, with
pairing $\pscal$, with norm $\|\cdot\|$, and with
duality mapping
$J = \partial \thalb\|\cdot\|^2$, where ``$\partial$'' stands
for the subdifferential operator from Convex Analysis.
Notation not explicitly defined here is standard and as in, e.g.,
\cite{Rocky,RockWets,Zalinescu}.

Recall that $A$ is a set-valued operator from $X$ to $X^*$, written
$A\colon X\To X^*$, if $A$ is a mapping
from $X$ to the power set of $X^*$, i.e.,
$(\forall x\in X)$ $Ax\subseteq X^*$. The \emph{graph} of $A$
is $\gra A = \menge{(x,x^*)\in X\times X^*}{x^*\in Ax}$. Such a mapping is
monotone if $(\forall (x,x^*)\in\gra A)$ $(\forall (y,y^*)\in\gra A)$
$\scal{x-y}{x^*-y^*}\geq 0$, and maximal monotone if it cannot be properly
extended without destroying monotonicity.
The \emph{domain} of $A$ is $\dom A = \menge{x\in X}{Ax\neq \emp}$
and the \emph{range} of $A$ is $\ran A = A(X) = \bigcup_{x\in X} Ax$.
The \emph{inverse} of $A$ is the operator $A^{-1}\colon X^*\To X$,
defined via $\gra A^{-1}=\menge{(x^*,x)\in X^*\times X}{(x,x^*)\in \gra
A}$. 

Monotone operators have turned out to be ubiquitous in modern
optimization and analysis;
see, e.g., \cite{Brezis,BuIu,RockWets,Simons08}. 
Due to their set-valuedness, there has
always been considerable interest to describe and study
monotone operators from a more classical point of view.
For ease of discussion, let us momentarily assume that $X$ is a Hilbert space.
A key tool is the so-called \emph{resolvent} $(A+\Id)^{-1}$ associated
with a given monotone operator $A$.
This resolvent is not only always
single-valued, but also \emph{firmly nonexpansive} (and thus
Lipschitz continuous); moreover, the resolvent has full domain $X$
precisely when $A$ is maximal monotone.
Resolvents can be used to parametrize the graph
of $A$, and the inverse-resolvent identity provides
a useful and elegant expression for the resolvent of $A^{-1}$
in terms of the resolvent for $A$.
More general resolvents have been studied.
Alber \cite{Alber96} and Kohsaka and Takahashi
\cite{KohTak04,KohTak08,KohTak08AM} 
initiated
the systematic study of 
resolvents based on the duality mapping $J$. Building on work by
Bregman \cite{Breg67} on generalized distances, Censor and Zenios 
analyzed proximal mappings \cite{CZ92} (see also \cite{CenZen}). 
For either generalization, it is known that every resolvent
is firmly nonexpansive.

\emph{The aim of this note is to present
a very general framework for resolvents and firmly nonexpansive
mappings in which the two classes coincide.} We also study
parametrizations of the graph, inverse resolvents,
and extensions of firmly nonexpansive mappings.
Various examples illustrate our results.

The paper is organized as follows. 
In Section~\ref{s:2}, we review the crucial characterization
due to Mart\'{\i}nez-Legaz (Fact~\ref{f:MLnew}) and 
then fix a monotone
operator $F$ upon which the various general notions are based.
Section~\ref{s:3} discusses
$F$-firmly nonexpansive mappings, and Section~\ref{s:4}
$F$-resolvents. It is then proved that
$F$-resolvents are $F$-firmly nonexpansive (Corollary~\ref{c:blabla});
the converse implication (Proposition~\ref{p:alaa}) is established in
Section~\ref{s:5}. The parametrization of the graph \`a la
Minty is obtained in Section~\ref{s:6}, while the resolvent
of the inverse is discussed in Section~\ref{s:7}.
Section~\ref{s:8} deals with the constructive extension
of a given $F$-firmly nonexpansive mapping.
The final Section~\ref{s:9} provides additional examples
and a foray into algorithms.

\section{Characterizations of maximality}
\label{s:2}

\begin{fact}[Mart\'{\i}nez-Legaz] \label{f:MLnew}
\emph{(See \cite[Theorem~8]{ML}.)}
Let $F\colon X\To X^*$ be a maximal monotone operator
such that its Fitzpatrick function \emph{\cite{FIT}}
\begin{equation} \label{e:Fitzfunc} 
X\times X^*\to\RX\colon
(x,x^*) \mapsto \sup_{(y,y*)\in\gra F} \big(\scal{x}{y^*} + \scal{y}{x^*}
-\scal{y}{y^*}\big)
\end{equation}
is real-valued,
and let $A\colon X\To X^*$ be monotone. Then the following hold.
\begin{enumerate}
\item \label{f:MLnewi}
If $A$ is maximal monotone, then $\ran(A+F)=X^*$.
\item \label{f:MLnewii}
If $F$ is single-valued, strictly monotone, and $\ran(A+F)=X^*$,
then $A$ is maximal monotone.
\end{enumerate}
\end{fact}

\begin{lemma} \label{l:3star}
Let $F\colon X\To X^*$ be a maximal monotone operator.
Then the Fitzpatrick function of $F$ is real-valued
$\Leftrightarrow$ $(\dom F)\times(\ran F) = X\times X^*$ and 
$F$ is 3*-monotone, i.e., 
\begin{equation}
\big(\forall (x,x^*) \in (\dom F)\times (\ran F)\big)\quad
\sup_{(y,y^*)\in\gra A}\scal{x-y}{y^*-x^*} <\pinf.
\end{equation}
\end{lemma}
\begin{proof}
``$\Rightarrow$'':
This follows from \cite[Corollary~3]{ML}. 
``$\Leftarrow$'': Clear.
\end{proof}

\begin{theorem} \label{t:MLnew}
Let $F\colon X\to X^*$ be maximal monotone, strictly monotone,
3*-monotone,
and surjective,
and let $A\colon X\To X^*$ be monotone.
Then
\begin{equation}
\text{$A$ is maximal monotone $\quad\Leftrightarrow\quad$
$\ran(A+\F)=X^*$.}
\end{equation}
\end{theorem}
\begin{proof}
Since $\dom F = X$, $\ran F = X^*$, and $F$ is 3*-monotone,
Lemma~\ref{l:3star} implies
that the Fitzpatrick function of $F$ is real-valued.
The characterization now follows from Fact~\ref{f:MLnew}.
\end{proof}

\begin{corollary} \label{c:ML}
Let $f\colon X\to\RR$ be G\^ateaux differentiable everywhere,
strictly convex, and cofinite,
and let $A\colon X\To X^*$ be monotone.
Then
$A$ is maximal monotone $\Leftrightarrow$
$\ran(A+\nabla f)=X^*$.
\end{corollary}
\begin{proof}
Indeed, $\dom \nabla f=X$ (by assumption),
$\nabla f$ is maximal monotone (as a subdifferential),
strictly monotone (as $f$ is strictly convex),
3*-monotone (as a subdifferential), and
$\ran \nabla f = \dom f^* = X^*$ (by assumption).
The result thus follows from Theorem~\ref{t:MLnew}.
\end{proof}

\begin{remark}
Some comments on Corollary~\ref{c:ML} are in order.
\begin{enumerate}
\item If $f$ is not cofinite, then the implication ``$\Rightarrow$''
fails: indeed, suppose that $X=\RR$, let
$f=\exp$, and set $A\equiv 0$.
Then $A$ is maximal monotone, yet $\ran(A+\nabla f) = \ran \nabla f
=\RPP \neq \RR$.
\item
If $f$ does not have full domain then the implication
``$\Leftarrow$'' fails:
this time, suppose that $X=\RR$, let $f$ be the negative
entropy function, and set $A = \Id|_{\RP}$.
Then $A+\nabla f = \Id +\nabla f $ is surjective
(which is seen either directly or from Corollary~\ref{c:Minty}),
but $A$ is not maximal monotone.
\end{enumerate}
\end{remark}

\begin{corollary}[Rockafellar] \label{c:Rocky}
\emph{(See  \cite[Corollary~on~page~78]{Rock70}, and also \cite{SimZal}
for another proof.)}
Suppose that $X$ is strictly convex and smooth,
and let $A\colon X\To X^*$ be monotone.
Then $A$ is maximal monotone $\Leftrightarrow$
$\ran(A+J)=X^*$.
\end{corollary}
\begin{proof}
To say that the Banach space $X$ is strictly convex and smooth
means precisely that
$\thalb\|\cdot\|^2$ is strictly convex and G\^ateaux differentiable.
Since $\thalb\|\cdot\|^2$ is cofinite
(the conjugate being the corresponding halved energy for the dual norm),
the result is clear from Corollary~\ref{c:ML}.
\end{proof}

Specializing Corollary~\ref{c:Rocky} further gives another
classical case.

\begin{corollary}[Minty] \emph{(See \cite{Minty}.)}
\label{c:Minty}
Suppose that $X$ is a Hilbert space, and
let $A\colon X\To X$ be monotone.
Then $A$ is maximal monotone $\Leftrightarrow$
$\ran(A+\Id)=X$.
\end{corollary}

\section{$F$-firmly nonexpansive operators}
\label{s:3}

From now on, we assume that
\boxedeqn{\label{e:standingassumption}
\text{
$F\colon X\to X^*$ is maximal monotone, strictly monotone,
3*-monotone, and surjective.}\\[+2 mm]
}
There are many examples of operators
satisfying our standing assumptions \eqref{e:standingassumption}
on $F$.

\begin{example} \label{ex:lots}
Each of the following describes a situation
where \eqref{e:standingassumption} holds.
\begin{enumerate}
\item \label{ex:lots:I}
$F=\Id$, when $X$ is a Hilbert space.
\item \label{ex:lots:J}
$F=J$, when $X$ is strictly convex and smooth. 
\item \label{ex:lots:p}
$F= \nabla \big(\tfrac{1}{p}\|\cdot\|^p\big)$, when 
$X$ is strictly convex and smooth, and $p\in\left]1,\pinf\right[$,
\item \label{ex:lots:f}
$F=\nabla f$, when 
$f\colon X\to\RR$ is differentiable everywhere, strictly convex,
and cofinite. 
\item \label{ex:lots:rot} $F$ is the counter-clockwise
rotator by an angle in $\left[0,\pi/2\right[$, when $X=\RR^2$. 
\end{enumerate}
\end{example}
\begin{proof}
It is clear that
\ref{ex:lots:I}--\ref{ex:lots:f} become increasingly
less restrictive;
for \ref{ex:lots:f}, 
the 3* monotonicity follows from \cite{BrezisHaraux} 
(see also \cite[Section~32.21]{Zeidler}). 
Finally, see \cite{BBW} for \ref{ex:lots:rot}. 
\end{proof}

\begin{definition}
Let $C\subseteq X$, and let
$T\colon C\to X$. Then $T$ is
\emph{$F$-firmly nonexpansive} if
\begin{equation}
(\forall x\in C)(\forall y\in C)\quad
\scal{Tx-Ty}{FTx-FTy} \leq \scal{Tx-Ty}{Fx-Fy}.
\end{equation}
\end{definition}

\begin{remark} While it is tempting to ponder
set-valued extension of $F$-firm nonexpansiveness, it turns
out that this leads one back to the single-valued case:
let  $T\colon X\To X$ satisfy
\begin{equation} \label{e:nosets}
\big(\forall (x,u)\in\gra T\big)\big(\forall (y,v)\in\gra T\big)\quad
\scal{u-v}{Fu-Fv} \leq \scal{u-v}{Fx-Fy},
\end{equation}
and suppose that $\{(x,u_1),(x,u_2)\}\subseteq \gra T$.
The monotonicity of $F$ and \eqref{e:nosets} yield
\begin{equation}
0\leq \scal{u_1-u_2}{Fu_1-Fu_2}\leq \scal{u_1-u_2}{Fx-Fx}=0.
\end{equation}
Hence $\scal{u_1-u_2}{Fu_2-Fu_2}=0$ and thus $u_1=u_2$ by
strict monotonicity of $F$.
\end{remark}

\begin{example}[classical firm nonexpansiveness]
\label{ex:classicalfn}
Suppose that $X$ is a Hilbert space
and that $F=\Id$. Let $C\subseteq X$ and let $T\colon
C\to X$. Then $T$ is $\Id$-firmly nonexpansive
$\Leftrightarrow$
\begin{equation}
(\forall x\in C)(\forall y\in C)\quad
\|Tx-Ty\|^2 \leq \scal{Tx-Ty}{x-y},
\end{equation}
i.e., $T$ is \emph{firmly nonexpansive} in the classical
Hilbert space sense (see, e.g.,  \cite{GoeKir,GoeRei}).
\end{example}

\begin{example}[``firmly nonexpansive type'']
Suppose that $X$ is strictly convex and smooth.
Let $C\subseteq X$ and let $T\colon C \to X$.
Following Kohsaka and Takahashi \cite{KohTak08}, we say
that the operator $T$ is of \emph{firmly nonexpansive type} if
\begin{equation}
(\forall x\in C)(\forall y\in C)\quad
\scal{Tx-Ty}{JTx-JTy} \leq \scal{Tx-Ty}{Jx-Jy}.
\end{equation}
\end{example}

\begin{example}[``$D$-firm'']
Let $f\colon X\to\RR$ be differentiable everywhere, strictly
convex, and cofinite,
let $C\subseteq X$, and let $T\colon C\to X$.
Following \cite{BBC03}, we say that the operator $T$
is \emph{$D$-firm} if
\begin{equation} \label{e:Dfirm}
(\forall x\in C)(\forall y\in C)\quad
\scal{Tx-Ty}{\nabla f(Tx)-\nabla f(Ty)} \leq \scal{Tx-Ty}{\nabla
f(x)-\nabla f(y)}.
\end{equation}
The ``$D$'' in $D$-firm stems from the fact that if we let
\begin{equation} \label{e:potential}
D\colon X\times X\to\RR\colon (x,y)\mapsto
f(x)-f(y)-\scal{x-y}{\nabla f(y)}
\end{equation}
be the \emph{Bregman distance} (see
\cite{Breg67,ButIus,CenZen} for further information) associated with $f$,
then \eqref{e:Dfirm} is equivalent to
\begin{equation}
(\forall x\in C)(\forall y\in C)
\quad
D(Tx,Ty) + D(Ty,Tx) \leq D(Tx,y)+D(Ty,x)-D(Tx,x)-D(Ty,y);
\end{equation}
see also \cite[Proposition~3.5(iv)]{BBC03}.
Note that if $X$ is strictly convex and smooth, and
$f = \tfrac{1}{2}\|\cdot\|^2$, then $T$
is $D$-firm $\Leftrightarrow$ $T$ is of firmly nonexpansive type.
In this sense, the notion of $D$-firmness is significantly more general
than that of firmly nonexpansive type.
\end{example}

In the next section, we turn to the construction of examples
of $F$-firmly nonexpansive operators.

\section{$F$-resolvents are $F$-firmly nonexpansive \ldots}
\label{s:4}

In the setting of Hilbert space, as in Example~\ref{ex:classicalfn},
it is well known that resolvents of monotone operators
are firmly nonexpansive.
More generally, operators that are of firmly nonexpansive type
or even $D$-firm may be obtained similarly. Most generally,
we will show in this section that $F$-resolvents give similarly rise
to $F$-firmly nonexpansive operators.

\begin{definition}
Let $A\colon X\To X^*$.
Then the composition
\begin{equation}
(A+F)^{-1}F
\end{equation}
is the \emph{$F$-resolvent} of $A$.
\end{definition}

\begin{proposition} \label{p:Fres}
Let $A\colon X\To X^*$,
let $T_A = (A+F)^{-1}F$ be its associated $F$-resolvent,
and let $x\in X$.
Then the following hold.
\begin{enumerate}
\item \label{p:Fresi}
$\dom T_A = F^{-1}(\ran(A+F))$ and $\ran T_A = \dom A$. 
\item \label{p:Fresiii}
$x\in T_Ax$ $\Leftrightarrow$ $0\in Ax$.
\item \label{p:Fresia}
If $A$ is monotone, then $T_A$ is at most single-valued
and $F$-firmly nonexpansive.
\item \label{p:Fresii}
If $A$ is monotone, then: $A$ is maximal monotone
$\Leftrightarrow$ $\dom T_A=X$.
\end{enumerate}
\end{proposition}
\begin{proof}
\ref{p:Fresi}:
$x\in\dom T_A$
$\Leftrightarrow$
$Fx\in\dom(A+F)^{-1}$
$\Leftrightarrow$
$Fx\in\ran(A+F)$
$\Leftrightarrow$
$x\in F^{-1}(\ran(A+F))$.
Furthermore, $\ran T_A = \dom F^{-1}(A+F) = \dom A$. 

\ref{p:Fresiii}:
$x\in T_Ax$
$\Leftrightarrow$
$Fx\in(A+F)x=Ax+Fx$
$\Leftrightarrow$
$0\in Ax$.

\ref{p:Fresia}: Suppose that $A$ is monotone.
Since $F$
is strictly monotone, it follows that $A+F$ is
strictly monotone, which in turn implies that $(A+F)^{-1}$
is at most single-valued. Since $F$ is single-valued,
we deduce that the composition $(A+F)^{-1}F$ is at most
single-valued. Using \ref{p:Fresi},
we set $C = \dom T_A = F^{-1}(\ran(A+F))$.
Let $y\in C$, i.e., $Fy\in\ran(A+F)$.
Then there exists $v\in X$ such that
$Fy \in (A+F)v$. Hence $Fy-Fv \in Av$ and
$v\in (A+F)^{-1}Fy=T_Ay$, i.e.,
$v=T_Ay$ and so
\begin{equation}
\label{e:kalamazoo}
(T_Ay,Fy-FT_Ay) \in \gra A.
\end{equation}
Let $z\in C$. A similar argugment shows that 
there exists $w = T_Az\in X$ such that
$Fz-Fw\in Az$ and $w=T_Az$. 
Since $A$ is monotone, 
$0 \leq \scal{v-w}{(Fy-Fv)-(Fz-Fw)}
= \scal{T_Ay-T_Az}{(Fy-Fz)-(FT_Ay-FT_Az)}$, 
i.e., 
\begin{equation}
\scal{T_Ay-T_Az}{FT_Ay-FT_Az} \leq
\scal{T_Ay-T_Az}{Fy-Fz}.
\end{equation}
This verifies that $T_A$ is $F$-firmly nonexpansive.

\ref{p:Fresii}:
Suppose that $A$ is monotone.
Using Theorem~\ref{t:MLnew},
the bijectivity of $F$, and \ref{p:Fresi},
we obtain
the equivalences:
$A$ is maximal monotone
$\Leftrightarrow$
$\ran(A+F)=X^*$
$\Leftrightarrow$
$F^{-1}(\ran(A+F))=X$
$\Leftrightarrow$
$\dom T_A = X$.
\end{proof}

\begin{corollary} \label{c:blabla}
Let $A\colon X\To X^*$ be maximal monotone, and
let $T_A = (A+F)^{-1}F$ be its associated $F$-resolvent.
Then $T_A\colon X\to X$ is $F$-firmly nonexpansive.
If $X$ is finite-dimensional, then $T_A$ is continuous.
\end{corollary}
\begin{proof}
In view of Proposition~\ref{p:Fres}, we only have to
establish the continuity of $T_A$ in the finite-dimensional case.
Since $F$ and $(A+F)^{-1}$ are single-valued maximal monotone
operators with full domain, it follows from
\cite[Theorem~12.63(c)]{RockWets} that they are continuous,
and so is their composition $(A+F)^{-1}F = T_A$. 
\end{proof}

\begin{example}
Let $F=\nabla f$ be as in Example~\ref{ex:lots}\ref{ex:lots:f}.
Then the $F$-resolvent of a maximal monotone operator $A$
becomes the ``$D$-resolvent'' considered in
\cite{Eck93,BBC03}, and the counterpart of
Proposition~\ref{p:Fres} is \cite[Proposition~3.8]{BBC03}.
If $A$ is a subdifferential operator, then one
obtains ``$D$-prox operators''; see, e.g., \cite{CZ92,BBC03}.
Finally, if $A=N_C$, where $C$ is a nonempty closed convex subset of $X$,
then we obtain Bregman projections; see, e.g.,
\cite{AlBut97}.
\end{example}

\begin{example}
Suppose that $X$ is strictly convex and smooth, and
let $F=J$ be as in Example~\ref{ex:lots}\ref{ex:lots:J}.
We then recover the resolvent $(A+J)^{-1}J$
(see, e.g., \cite{Kassay85,KohTak08}), and the counterpart
of Proposition~\ref{p:Fres} is \cite[Lemma~2.3]{KohTak08}.
If $A$ is specialized to the normal cone operator $N_C$,
where $C$ is a nonempty closed convex subset of $X$,
then the resolvent becomes the generalized projection operators
studied, e.g., in \cite{Alber96,KohTak08}.
\end{example}

\begin{example}[Minty-Rockafellar] \label{ex:MintRock}
Suppose $X$ is a Hilbert space and $A$ is maximal
monotone. Then the standard resolvent $(A+\Id)^{-1}$ is firmly
nonexpansive and it has full domain.
This is classical and goes back to Minty \cite{Minty}
and to Rockafellar \cite{Rocky76}.
\end{example}

\section{\ldots\ and vice versa}
\label{s:5}

Eckstein and Bertsekas \cite{EckBer} observed that the converse
of Example~\ref{ex:MintRock} holds, i.e.,
that every firmly nonexpansive operator (with full domain)
must be the resolvent of the corresponding
(maximal) monotone operator.
As we now show, this is also the case for $F$-firmly nonexpansive
operators.

\begin{proposition} \label{p:alaa}
Let $C\subseteq X$, let $T\colon C\to X$, and
set $A_T = FT^{-1} - F$.
Then the following hold.
\begin{enumerate}
\item \label{p:alaa:i}
The $F$-resolvent of $A_T$ is $T$.
\item \label{p:alaa:ii}
If $T$ is $F$-firmly nonexpansive,
then $A_T$ is monotone.
\item \label{p:alaa:iii}
If $T$ is $F$-firmly nonexpansive,
then: $C=X$ $\Leftrightarrow$ $A_T$ is maximal monotone.
\end{enumerate}
\end{proposition}
\begin{proof}
\ref{p:alaa:i}:
$A_T = F T^{-1}-F$
$\Rightarrow$ $A_T+F = FT^{-1}$
$\Rightarrow$ $(A_T+F)^{-1} = (F T^{-1})^{-1} = TF^{-1}$
$\Rightarrow$ $(A_T+F)^{-1}F = TF^{-1}F = T$.

\ref{p:alaa:ii}:
Suppose that $T$ is $F$-firmly nonexpansive.
Take $(u,u^*),(v,v^*)$ in $\gra A_T$.
Then $u^*\in A_Tu = FT^{-1}u-Fu$
$\Leftrightarrow$
$u^*+Fu \in FT^{-1}u$
$\Leftrightarrow$
$u\in (FT^{-1})^{-1}(u^*+Fu)$
$\Leftrightarrow$
$u =  TF^{-1}(u^*+Fu)$, and analogously $v=TF^{-1}(v^*+Fv)$.
Since $T$ is $F$-firmly nonexpansive, we estimate
\begin{align}
\scal{u-v}{Fu-Fv} &=
\scal{TF^{-1}(u^*+Fu)-TF^{-1}(v^*+Fv)}{FTF^{-1}(u^*+Fu)-FTF^{-1}(v^*+Fv)}\notag\\
&\leq
\scal{TF^{-1}(u^*+Fu)-TF^{-1}(v^*+Fv)}{FF^{-1}(u^*+Fu)-FF^{-1}(v^*+Fv)}
\notag\\
&= \scal{u-v}{(u^*+Fu)-(v^*+Fv)}.
\end{align}
Hence, $0\leq \scal{u-v}{u^*-v^*}$, as required.

\ref{p:alaa:iii}:
Suppose that $T$ is $F$-firmly nonexpansive.
By \ref{p:alaa:ii}, $A_T$ is monotone.
Using \ref{p:alaa:i} and Proposition~\ref{p:Fres}\ref{p:Fresii},
we obtain: $A_T$ is maximal monotone $\Leftrightarrow$
$\dom T= C =X$.
\end{proof}

\begin{corollary} \label{c:plane}
Let $A\colon X\To X^*$ with associated $F$-resolvent $T_A = (A+F)^{-1}F$,
let $C\subseteq X$, let $T\colon C\to X$, and set $A_T = FT^{-1}-F$.
Assume that $T_A = T$; equivalently, that $A_T = A$.
Then $A$ is (maximal) monotone
$\Leftrightarrow$
$T$ is $F$-firmly nonexpansive (and $C=X$).
\end{corollary}
\begin{proof}
Combine Proposition~\ref{p:Fres}
and Proposition~\ref{p:alaa}.
\end{proof}

Specializing to $F=J$, where $X$ is strictly convex and smooth,
one obtains the following result related
to \cite[Proposition~3.1]{KohTak08AM}.

\begin{corollary}[Kohsaka-Takahashi]
Suppose that $X$ is strictly convex and smooth
and that $F=J$. 
Let $C\subseteq X$, let $T\colon C\to X$,
and set $A_T = JT^{-1}-J$.
Then $T$ is $J$-firmly nonexpansive
$\Leftrightarrow$
$A_T$ is monotone.
\end{corollary}

In the setting
of Hilbert space, Corollary~\ref{c:plane} recovers
the following result, which appeared
first in \cite[Theorem~2]{EckBer}. 

\begin{corollary}[Eckstein-Bertsekas]
Suppose that $X$ is a Hilbert space and
that $F=\Id$,
let $A\colon X\To X$, and
denote the $\Id$-resolvent of $A$ by $T_A$.
Then $A$ is (maximal) monotone
$\Leftrightarrow$
$T_A$ is firmly nonexpansive (with full domain). 
\end{corollary}

\section{Minty parametrization}
\label{s:6}

\begin{theorem}[$F$-Minty parametrization]
\label{t:mintyparmap}
Let $A\colon X\To X^*$ be monotone,
let $T_A = (A+F)^{-1}F$ be its associated $F$-resolvent,
and set $C=\dom T_A$. Then
\begin{equation} 
\label{e:mintyparmap}
\Psi\colon C \to \gra A \colon x \mapsto
\big(T_Ax,Fx-FT_Ax\big)
\end{equation}
is a bijection with 
\begin{equation}
\label{e:mintyparmapinv}
\Psi^{-1}\colon \gra A\to C\colon
(u,u^*)\mapsto F^{-1}(u^*+Fu).
\end{equation}
Moreover, the following hold.
\begin{enumerate}
\item \label{t:mintyparmap:i}
If $F,F^{-1},T_A$ are continuous, then so are $\Psi$ and $\Psi^{-1}$. 
\item \label{t:mintyparmap:ii}
If $X$ is finite-dimensional and $A$ is maximal monotone, 
then $F,F^{-1},T_A,\Psi,\Psi^{-1}$ are continuous.
\item \label{t:mintyparmap:iii}
If $X$ is finite-dimensional and $F$ is linear, then
$F,F^{-1},T_A,\Psi,\Psi^{-1}$ are Lipschitz continuous.
\end{enumerate}
\end{theorem}
\begin{proof}
It follows from \eqref{e:kalamazoo} that 
$(\forall y\in C)$ $(T_Ay,Fy-FT_Ay)\in\gra A$.
Hence $\ran\Psi \subseteq\gra A$.
Now take $(u,u^*)\in\gra A$
and set $x = F^{-1}(u^*+Fu)$.
Then $u^*\in Au$ $\Rightarrow$ $u^*+Fu\in (A+F)u$
$\Rightarrow$ $x=F^{-1}(u^*+Fu)\in F^{-1}(\ran(A+F)) = \dom T_A = C$
by Proposition~\ref{p:Fres}\ref{p:Fresi}.
Furthermore, $T_Ax=(A+F)^{-1}FF^{-1}(u^*+Fu) =
(A+F)^{-1}(u^*+Fu)=u$ and thus
$Fx-FT_Ax = FF^{-1}(u^*+Fu)-Fu = u^*$.
Therefore, $(u,u^*)=\Psi(x)$ and hence $\ran\Psi=\gra A$.
On the other hand, let $y$ and $z$ be in $C$
such that $\Psi(y)=\Psi(z)$.
Then $T_Ax=T_Ay$ and $Fx-FT_Ax=Fy-FT_Ay$, hence
that $Fx=Fy$ and thus $x=y$.
It follows that $\Psi$ is injective.
Altogether, $\Psi$ is a bijection between $C$ and $\gra A$.
The beginning of this proof implies the formula for $\Psi^{-1}$.
We now turn to the continuity assertions. 

\ref{t:mintyparmap:i}: This statement is clear from
the formulae \eqref{e:mintyparmap} and \eqref{e:mintyparmapinv}. 

\ref{t:mintyparmap:ii}: 
Suppose that $X$ is finite-dimensional.
In Corollary~\ref{c:blabla}, we
observed that $T_A$ is continuous; 
by using once again
\cite[Theorem~12.63(c)]{RockWets}, we obtain continuity
of $F$ and $F^{-1}$. Now apply \ref{t:mintyparmap:i}. 

\ref{t:mintyparmap:iii}: 
It is clear that $F$ and $F^{-1}$ are Lipschitz continuous. 
Denote the smallest eigenvalue of the symmetric part of $F$ by
$\lambda$. Then $\lambda>0$. Since $T_A$ is $F$-firmly nonexpansive
by Proposition~\ref{p:Fres}\ref{p:Fresia}, we estimate
\begin{align*}
(\forall x\in C)(\forall y\in C)\quad \lambda\|T_Ax-T_Ay\|^2 &\leq \scal{T_Ax-T_Ay}{F(T_Ax-T_Ay)} \\
&=\scal{T_Ax-T_Ay}{FT_Ax-FT_Ay} \\
&\leq \scal{T_Ax-T_Ay}{Fx-Fy}\\
&\leq \|T_Ax-T_Ay\|\,\|F\|\,\|x-y\|;
\end{align*}
consequently,
$\|T_Ax-T_Ay\|\leq (\|F\|/\lambda)\|x-y\|$. 
The formulae \eqref{e:mintyparmap} and \eqref{e:mintyparmapinv} show
that $\Psi$ and $\Psi^{-1}$ are Lipschitz continuous as well. 
\end{proof}

\begin{remark}
When $F=J$, 
the inclusion $\ran\Psi \subseteq \gra A$ in
Theorem~\ref{t:mintyparmap} was already noted
by Kohsaka and Takahashi (see \cite[page~242]{KohTak04}).
\end{remark}

\section{Resolvent of the inverse}
\label{s:7}

In this section, we discuss the possibility of
computing the $F^{-1}$-resolvent of $A^{-1}$
in terms of the $F$-resolvent of $A$.

\begin{theorem}[inverse-resolvent fixed point equation]
Let $A\colon X\To X^*$ be monotone,
let $T_A=(A+F)^{-1}F$ be the its associated $F$-resolvent,
let $T_{A^{-1}} = (A^{-1}+F^{-1})^{-1}F^{-1}$
be the $F^{-1}$-resolvent of $A^{-1}$, let
$x^*\in\dom T_{A^{-1}} = F(\ran(A^{-1}+F^{-1})$,
and let $y^*\in X^*$.
Then 
\begin{equation} \label{e:irife}
y^*=T_{A^{-1}}x^* \;\;\Leftrightarrow\;\; y^* = F\bigg(F^{-1}x^*-T_A\Big(F^{-1}\big(y^*+F(F^{-1}x^*-F^{-1}y^*)\big) \Big)
\bigg).
\end{equation}
\end{theorem}
\begin{proof}
The identity for $\dom T_{A^{-1}}$ follows from
Proposition~\ref{p:Fres}\ref{p:Fresi}.
For convenience, set $x=F^{-1}x^*$ and $y=F^{-1}y^*$.
We then have the equivalences
{\allowdisplaybreaks
\begin{align}
y^* = T_{A^{-1}}x^* &
\Leftrightarrow y^* = (F^{-1}+A^{-1})^{-1}F^{-1}x^*\notag\\
&\Leftrightarrow F^{-1}x^* \in (F^{-1}+A^{-1})y^*\notag\\
&\Leftrightarrow x-y \in A^{-1}y^*\notag\\
&\Leftrightarrow y^*\in A(x-y)\notag\\
&\Leftrightarrow y^*+F(x-y)\in (A+F)(x-y)\notag\\
&\Leftrightarrow x-y\in (A+F)^{-1}FF^{-1}\big(y^*+F(x-y)\big)\notag\\
&\Leftrightarrow x-y= T_AF^{-1}\big(y^*+F(x-y)\big)\notag\\
&\Leftrightarrow y= x-T_AF^{-1}\big(y^*+F(x-y)\big)\notag\\
&\Leftrightarrow y^*= F\Big(x-T_AF^{-1}\big(y^*+F(x-y)\big)\Big),
\end{align}
}
and this last identity is in turn equivalent to 
the right side of \eqref{e:irife}.
\end{proof}

\begin{corollary} Suppose that $F$ is linear, 
let $A\colon X\To X^*$ be monotone,
let $T_A = (A+F)^{-1}F$ be its associated $F$-resolvent,
and let $T_{A^{-1}} = (A^{-1}+F^{-1})^{-1}F^{-1}$ 
be the $F^{-1}$-resolvent of $A^{-1}$. 
Then
\begin{equation}
T_{A^{-1}} = \Id - FT_AF^{-1}.
\end{equation}
\end{corollary}

In the classical Hilbert space setting of
Example~\ref{ex:MintRock}, one recovers
the following well known result \cite[Lemma~12.14]{RockWets}.

\begin{corollary}[inverse-resolvent identity]
Suppose that $X$ is a Hilbert space and that $F=\Id$.
Let $A\colon X\To X^*$ be maximal monotone.
Then $T_{A^{-1}} = \Id - T_A$, i.e., 
\begin{equation}
(A^{-1}+\Id)^{-1} = \Id - (A+\Id)^{-1}.
\end{equation}
\end{corollary}

\section{Constructive extension}
\label{s:8}

We now describe how $F$-firmly nonexpansive operators
can be extended to the whole space. This technique
was recently utilized in \cite{BaWa08} in the setting
of Hilbert spaces. 

\begin{theorem} \label{t:kalala}
Let $C\subseteq X$, and let $T\colon C\to X$
be $F$-firmly nonexpansive. Proceed as follows.
\begin{dingautolist}{192}
\item Set $A = FT^{-1}-F$.
\item Denote the Fitzpatrick function of $A$ (see \eqref{e:Fitzfunc}) 
by $\Phi$.
\item Compute 
\begin{equation}
\Psi\colon (x,x^*)\mapsto \min_{(y+z,y^*+z^*) = 2(x,x^*)}
\big( \thalb\Phi(y,y^*) + \thalb\Phi^*(z^*,z) + \tfrac{1}{8}\big(\|y-z\|^2
+ \|y^*-z^*\|^2\big)\big),
\end{equation}
which is the {proximal average} \emph{\cite{BGLW}}
between $\Phi$ and $\Phi^*$ (with the variables transposed). 
\item Define $\widetilde{A}\colon X\To X^*$ via
\begin{equation}
\gra \widetilde{A} = \menge{(x,x^*)\in X\times
X^*}{\Psi(x,x^*)=\scal{x}{x^*}}.
\end{equation}
\item Set $\widetilde{T} = (\widetilde{A}+F)^{-1}F$.
\end{dingautolist}
Then $\widetilde{T}\colon X\to X$ is $F$-firmly nonexpansive
and it extends $T$ to the entirety of $X$. 
\end{theorem}
\begin{proof}
By Proposition~\ref{p:alaa}, $A$ is monotone. 
Hence, using \cite[Fact~5.6 and Theorem~5.7]{BaWapress},
we see that $\widetilde{A}$ is a maximal monotone extension
of $A$. Theorem~\ref{t:mintyparmap} and Proposition~\ref{p:Fres} 
now show that $\widetilde{T}$ is an $F$-firmly nonexpansive 
extension of $T$ to the entire space $X$. 
\end{proof}

\begin{remark}
Let us comment on Theorem~\ref{t:kalala} further when
$X$ is a real Hilbert space. 
\begin{enumerate}
\item In this case, Theorem~\ref{t:kalala}
becomes \cite[Theorem~3.1]{BaWa08}. 
\item As explained in \cite[Theorem~3.6]{BaWa08}, one may
use Theorem~\ref{t:kalala} to obtain a \emph{constructive}
Kirszbraun-Valentine extension of a given nonexpansive operator.
\end{enumerate}
\end{remark}

\section{Examples}
\label{s:9}

We begin with the $F$-resolvent of the identity,
where $F$ is a counter-clockwise rotator in the Euclidean plane.

\begin{example} \label{ex:roti}
Suppose that $X = \RR^2$, 
let $\theta \in\left[0,\tfrac{\pi}{2}\right[$, and set
\begin{align}
F &=\begin{pmatrix}
\cos\theta & -\sin\theta \\
\sin \theta  & \cos\theta
\end{pmatrix}.
\end{align}
Then 
\begin{equation} 
(\Id+F)^{-1}F=\thalb\begin{pmatrix}
1 & -\tfrac{\sin\theta}{1+\cos\theta} \\[+2 mm]
\tfrac{\sin\theta}{1+\cos\theta}  & 1
\end{pmatrix}.\end{equation}
\end{example}

The most important example of a standard resolvent is the projector
onto a nonempty closed convex set $C$, which arises
as the resolvent of the normal cone operator
$N_C = \partial\iota_C$. As it turns out, 
a generalized projector is obtained in the general 
$F$-resolvent setting. 

\begin{theorem}[$F$-projector]
\label{t:superP_C}
Let $C \subseteq X$ be nonempty, closed, and convex,
denote the $F$-resolvent of $N_C$ by $P_C$,
and assume that $y\in\inte C$. 
Then $\ran P_C = \Fix P_C = C$, $P_C^2 = P_C$, and 
$P_C^{-1}y = \{y\}$. 
\end{theorem}
\begin{proof}
Note that $\ran P_C = \dom N_C = C$ by
Proposition~\ref{p:Fres}\ref{p:Fresi};
furthermore,
$\Fix P_C = N_C^{-1}0 = C$ by Proposition~\ref{p:Fres}\ref{p:Fresiii}. 
Finally, since $y\in\inte C$, $N_Cy=\{0\}$ and therefore
$P_C^{-1}y = \big((N_C+F)^{-1}F\big)^{-1}y=F^{-1}(N_C+F)y = 
F^{-1}(0+Fy)=y$.
\end{proof}

For the purpose of illustration, 
let us now compute some generalized projectors when
$F$ is the rotator from Example~\ref{ex:roti}.
The following result is clear from Theorem~\ref{t:superP_C}. 

\begin{example}
Suppose that $X=\RR^2$, let $\theta$ and $F$ be
as in Example~\ref{ex:roti},
let $C\subseteq \RR^2$, and denote the $F$-resolvent of $N_C$
by $P_C$. \begin{enumerate}
\item If $C=\{0\}$, then $P_C=0$.
\item If $C=\RR^2$, then $P_C=\Id$.
\end{enumerate}
\end{example}

\begin{example}
Suppose that $X=\RR^2$, let $\theta$ and $F$ be
as in Example~\ref{ex:roti},
set $C=\RR\times\{0\}$, and denote the $F$-resolvent of $N_C$
by $P_C$. Then 
\begin{equation} 
P_C =\begin{pmatrix}
1 & -\tan\theta \\
0 & 0
\end{pmatrix}.\label{Ex1}\end{equation}
\end{example}
\begin{proof}
Let $x=(x_1,x_2) \in\RR^2$, and set $y=P_Cx$. Then
$y\in C$ and $Fx\in N_Cy+Fy=(\{0\}\times \RR)+Fy$.
Thus $F(x-y)\in\{0\}\times \RR.$
Write $y=(y_1,0)$. We then have
$(x_1-y_1)\cos{\theta}-x_2\sin{\theta}=0$.
Hence $y_1=x_1-x_2\tan{\theta}$, and \eqref{Ex1} holds.
\end{proof}

\begin{example}
Suppose that $X=\RR^2$, let $\theta$ and $F$ be
as in Example~\ref{ex:roti},
let $C=\menge{x\in\RR^2}{\|x\|\leq 1}$ be the closed unit ball, 
denote the $F$-resolvent of $N_C$ by $P_C$, and
set $\alpha = \sqrt{\|x\|^2-\sin^2\theta}-\cos\theta$. Then 
\begin{equation}
P_C\colon \RR^2 \to C \colon
x\mapsto \begin{cases}
 x, &\text{if $x\in C$;}\\
\tfrac{1}{\|x\|^2}(\Id+\alpha F) x,\label{B:2}
&\text{if $x\notin C$.}
\end{cases}
\end{equation}
Moreover, 
\begin{equation} 
\big(\forall z\in\RR^2\big)\quad 
P_C^{-1}z
=\begin{cases}z,&\text{if $\|z\|<1$;}\\
z+\RP\cdot F^*z,&\text{if $\|z\|=1$;}\\
\emp,&\text{otherwise}.\end{cases}
\label{B:7}\end{equation}
\end{example}
\begin{proof} 
Let $x\in\RR^2$. We consider two cases.\\
\emph{Case 1}: $\|x\|\leq 1$.
Then $x\in C$ and so $P_Cx=x$ by Theorem~\ref{t:superP_C}.\\
\emph{Case 2}: 
$\|x\|>1$. Set $y=P_Cx=(N_C+F)^{-1}Fx$.
Assume that $\|y\|<1$. Then 
$Fx=N_Cy+Fy=0+Fy=Fy$. Hence $x=y$, which is absurd. Thus
\begin{equation}
\|y\|=1.
\end{equation}
and therefore
\begin{equation} 
\label{B:1}
N_C(y)=\RP \cdot y.
\end{equation}
It follows that there exists $\alpha\in\RP$ 
such that $Fx=\alpha y+Fy=(\alpha \Id+F)y$.
Since $x\neq y$, we see that $\alpha>0$.
Moreover, by the orthogonality of $F$, it follows that 
\begin{equation}
(\alpha\Id + F)^{-1}F = 
\tfrac{1}{\alpha^2+2\alpha\cos\theta+1}(\alpha F+\Id)
\end{equation}
and that
\begin{equation}
F^*(\alpha\Id+F^*)^{-1}(\alpha\Id+F)^{-1}F = 
\tfrac{1}{\alpha^2+2\alpha\cos\theta+1}\Id.\label{Br:1}
\end{equation}
Since $1=\|y\|=\|(\alpha\Id+F)^{-1}Fx\|$, we thus have 
$\alpha^2 + 2\alpha\cos\theta + 1 = \|x\|^2$ and
hence 
$\alpha = -\cos\theta + \sqrt{\cos^2\theta +\|x\|^2-1}
=-\cos\theta + \sqrt{\|x\|^2-\sin^2\theta}$.
Consequently,
\begin{equation}
y = (\alpha\Id+F)^{-1}Fx = \tfrac{1}{\|x\|^2}(\alpha F+\Id)x,
\end{equation}
which yields \eqref{B:2}.
Now let $z\in\RR^2$.
In view of Theorem~\ref{t:superP_C}, it suffices
to consider the case when $z\in\bd C$, i.e., 
$\|z\|=1$ and thus $N_Cz = \RP\cdot z$.
Then $P_C^{-1}z = F^{-1}(N_C+F)z = z + F^*\RP\cdot z
= z + \RP\cdot F^*z$. 
\end{proof}

\begin{example} \label{ex:Sunday}
Suppose that $X$ is a Hilbert space and
that $F=\nabla \tfrac{1}{p}\|\cdot\|^p$,
where $p\in\left]1,\pinf\right[$.
Let $x\in X$ and set
\begin{equation}
k_p(x) = 
\begin{cases}
0, &\text{if $x=0$;}\\
\text{the unique solution of $k^{p-1} + {k}/{\|x\|^{p-2}} = 1$ in
$\left]0,1\right[$}, &\text{if $x\neq 0$.}
\end{cases}
\end{equation}
Let $T_p=(\Id+F)^{-1}F$ be the $F$-resolvent of $\Id$.
Then $T_p(x)=k_p(x)x$.
Moreover, 
\begin{equation}
\lim_{p\to 1^+} T_p(x) = 0 \quad\text{and}\quad
\lim_{p\to\pinf} T_p(x) = \begin{cases}
0, &\text{if $\|x\|<1$;}\\
x, &\text{if $\|x\|\geq 1$.}
\end{cases}
\end{equation}
\end{example}
\begin{proof}
The statements are clear if $x=0$, so we assume that $x\neq 0$.
Set $y=T_p(x)$. Then
$y\neq 0$, $Fx=\|x\|^{p-2}x$ and $Fy=\|y\|^{p-2}y$.
Furthermore, 
$Fx\in(\Id+F)y=y+Fy$ 
$\Leftrightarrow$
$\|x\|^{p-2}x = (1+\|y\|^{p-2})y$,
which implies that $y=k x$, where $k\in\RPP$
satisfies $\|x\|^{p-2} = k+ k^{p-1}\|x\|^{p-2}$
$\Leftrightarrow$
$k^{p-1}+{k}/{\|x\|^{p-2}} = 1$.
The remaining statements follow using Calculus.
\end{proof}

\begin{remark}
Consider Example~\ref{ex:Sunday} when $X$ is finite-dimensional.
By Corollary~\ref{c:blabla}, $T_p$ is continuous; 
however, the limiting (in the pointwise sense) operator 
$\lim_{p\to\pinf} T_p$ is not continuous.
\end{remark}

We now turn to an algorithmic result
on iterating $F$-resolvents. 

\begin{theorem} \label{t:Banach}
Suppose that $X$ is a Hilbert space and that $F$ is linear. 
Let $T=(\Id+F)^{-1}F$ be the $F$-resolvent of $\Id$, 
let $x_0\in X$, and set 
$(\forall\nnn)$ $x_{n+1}=Tx_n$. 
Then $\|T\|<1$ and hence $x_n\to 0$.
\end{theorem}
\begin{proof}
Set $\alpha = 1/\|F\|$. Then $(\forall x\in X)$ $\|Fx\|\leq \|x\|/\alpha$;
equivalently, 
\begin{equation} \label{e:rafal} 
(\forall y\in X)\quad
\|F^{-1}y\| \geq \alpha\|y\|. 
\end{equation}
Observe that 
\begin{equation}
T = (\Id+F)^{-1}F = (\Id + F^{-1})^{-1}, 
\end{equation}
let $x\in X$, and set $y=Tx$.
Then $x=T^{-1}y = (\Id+F^{-1})y$; thus,
by monotonicity of $F^{-1}$ and \eqref{e:rafal}, we obtain 
\begin{equation}
\|x\|^2 = \|y+F^{-1}y\|^2
= \|y\|^2 + \|F^{-1}y\|^2 + 2\scal{y}{F^{-1}y}
\geq \|y\|^2 + \|F^{-1}y\|^2
\geq (1+\alpha^2)\|y\|^2.
\end{equation}
Hence $\|y\|^2=\|Tx\|^2 \leq \|x\|^2/(1+\alpha^2)$ and so
\begin{equation}
\|T\| \leq 1/\sqrt{1+\alpha^2} < 1. 
\end{equation}
By the Banach Contraction Mapping Principle, 
we see that $(x_n)_\nnn = (T^nx_0)_\nnn$ converges in norm
to $0$, which is the unique fixed point of $T$. 
\end{proof}

\begin{remark}
Let us conclude by interpreting Theorem~\ref{t:Banach}
and outlining possible future research directions. 
Resolvent iterations are important for finding zeros
of subdifferential operators --- that is, minimizers --- or
more generally for finding zeros of maximal monotone operators. 
When $F=\Id$, this brings us to the  classical setting
of the proximal point algorithm \cite{Mart70,Rocky76}; 
when $F=J$, where $X$ is uniformly convex and uniformly smooth,
see \cite{KohTak08} and references therein, and
when $F=\nabla f$, this goes back to \cite{CZ92}.
It would be very interesting to build a general
convergence theory for iterating $F$-resolvents.
The difficulty lies in the absence of a potential
function like the Bregman distance \eqref{e:potential}.
However, Theorem~\ref{t:Banach} shows that
it may be possible to create a theory in the present
general framework, since this result shows that
resolvent iterations do converge to the unique zero
of the maximal monotone operator $\Id$. 
This promises to be an exciting topic for
further research. 
\end{remark}

\section*{Acknowledgment}

Heinz Bauschke was partially supported by the Natural Sciences and
Engineering Research Council of Canada and
by the Canada Research Chair Program.
Xianfu Wang was partially supported by the Natural
Sciences and Engineering Research Council of Canada.

\small 



\end{document}